\documentstyle[12pt,amssymb,amsmath,epsf,graphics,color]{article}

\def\X{{\cal S}}

\def\R {{\Bbb R}}
\def\N {{\Bbb N}}
\def \Z {{\Bbb Z}}
\def \S {{\Bbb S}}
\def\L {{\Bbb L}}
\newcommand{\D}{{\mathbb{D}}}

\def\cqfd{\hfill$\Box$}
\newcommand{\boM}{{\mathcal{M}}}
\newcommand{\boA}{{\mathcal{A}}}
\newcommand{\boB}{{\mathcal{B}}}
\newcommand{\boS}{{\mathcal{S}}}

\newcommand{\Ome}{\Omega}
\newtheorem{theorem}{Theorem}
\newtheorem{proposition}{Proposition}
\newtheorem{definition}{Definition}
\newtheorem{lemma}{Lemma}
\newtheorem{corollary}{Corollary}
\newtheorem{remark}{Remark}
\newtheorem{claim}{Claim}

\newcommand{\myskip}[1]{}

\begin{document}
\title{A quasi-periodic minimal surface}
\author{Laurent Mazet \\Martin Traizet}

\maketitle

{\em Abstract : we construct a properly embedded
minimal surface in the flat product $\R^2\times\S^1$
which is quasi-periodic but is not periodic.}

\section{Introduction}\label{secintro}

Quasi-periodicity is a popular subject in both mathematics and physics.
Probably the most famous examples are Penrose's quasi-periodic tiling, 
and quasi-periodic cristals.

\medskip

Recall that a planar tiling ${\cal T}$
is quasi-periodic if any finite  
part of the tiling repeats infinitely many often. In other words, for arbitrary
$R>0$, the tiling ${\cal T}$ countains an infinite number of
translation copies of ${\cal T}\cap B(0,R)$ where $B(0,R)$ is the ball
of radius $R$ centered at $0$.

\medskip

Of course, for minimal surfaces, it is too much to ask that a part of
the surface repeats exactly, because then by analytic continuation
the whole surface would be periodic.
We are thus led to the following definition, which was
suggested to the authors by H. Rosenberg.
\begin{definition}
A complete minimal surface $M$ in $\R^3$ is quasi-periodic if
there exists a diverging sequence of translations
$(T_n)_{n\in\N}$,
such that $T_n(M)$ converges smoothly to $M$ on compact
subsets of $\R^3$.
\end{definition}
While writing this paper, the authors discovered that the same notion
had been introduced by Meeks, Perez and Ros in a recent paper
\cite{meeks-perez-ros}, altough they call it translation-periodic.

\medskip

Of course a periodic minimal surface is quasi-periodic. A natural, and open,
question is whether there exists quasi-periodic minimal surfaces which
are not periodic.
In this paper we answer this question when the ambiant space is
the flat product $\R^2\times\S^1$ instead of $\R^3$.
The definition of quasi-periodicity is exactly the same in this case.
\begin{theorem}
There exists a complete embedded minimal surface in $\R^2\times\S^1$
which is quasi-periodic but is not periodic.
This surface has bounded curvature, infinite total curvature,
infinite genus, infinitely many ends and two limit ends.
\end{theorem}
Let us now explain informally how this surface is constructed.
H. Karcher has constructed a family of doubly periodic minimal surfaces
in $\R^3$
which he called the ``toroidal halfplane layers''
\cite{ka4}.
They were the first complete, properly embedded, doubly periodic minimal
surfaces to be found since H. Scherk's classical example.
The toroidal halfplane layers have two
periods : a horizontal period $T$ and a vertical period $(0,0,1)$.
We may identify the quotient of $\R^3$ by the vertical period $(0,0,1)$ with
$\R^2\times\S^1$. So the toroidal halfplane layers project to simply
periodic minimal surfaces in $\R^2\times\S^1$, with period $T$.
They have genus zero.

\medskip

A very succesful heuristic to construct new examples of minimal
surfaces is to start from a simple example, and to complicate it
by adding handles. One can start from a very symmetric example
and break the symmetries by adding handles at suitable places.

\medskip

Several people have added handles to H. Karcher's toroidal halfplane
layers.
The first one was F. Wei \cite{wei2}. He was able to add one handle
per fundamental piece in a periodic way. The
resulting surfaces have infinite genus in $\R^2\times\S^1$ and are
periodic.
W. Rossman, E. Thayer and M. Wolgemuth \cite{rossman-thayer-wohlgemuth}
have added handles in various
ways to the toroidal halfplane layers, still requiring periodicity.
Recently, the first author was able to add one single handle to the
toroidal halfplane layers, without requiring horizontal periodicity.
The resulting surfaces in $\R^2\times\S^1$ have genus one and are
not periodic anymore.

\medskip

In this paper we add an infinite number of handles to the
toroidal halfplane layers, so the resulting surface have infinite genus,
but without requiring horizontal periodicity. In fact the placement of the
handles will be prescribed by a sequence of integers $(p_i)_{i\in\Z}$.
If this sequence if quasi-periodic but not periodic, the resulting surface
will be quasi-periodic but not periodic.

\medskip

\begin{figure}
\begin{center}
\epsfysize=8.5cm
\epsffile{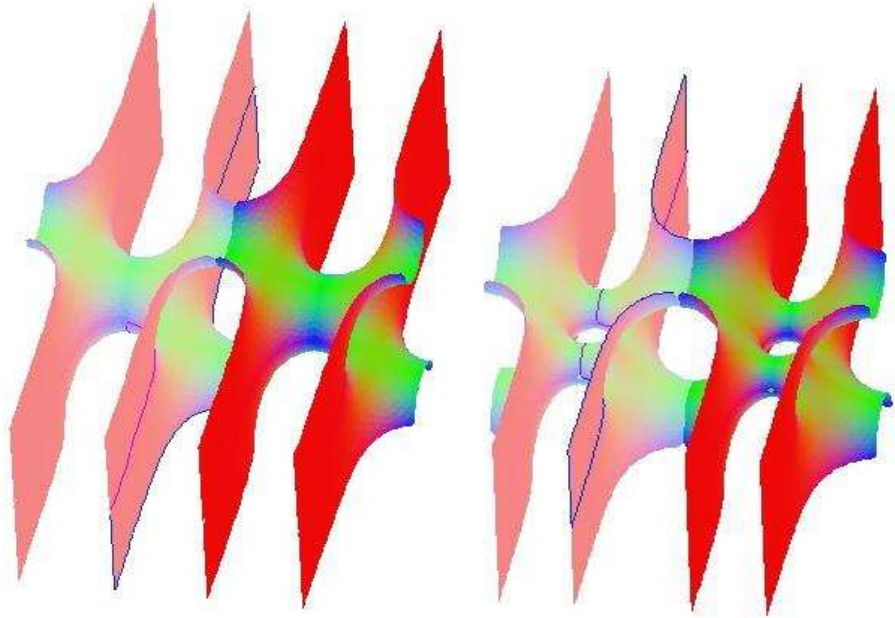}
\end{center}
\caption{Left : one of the Karcher's toroidal halfplane layers. Right
  : one of Wei surfaces. A fundamental domain is highlighted for
  each. Both surfaces extend periodically vertically and
  horizontally. The fundamental domains of these two surfaces are the
  basic building blocks for the surface we construct: we assemble
  them in a quasi periodic, non periodic way. Computer images made by
  the authors using J. Hoffman's MESH software.} 
\end{figure}

To construct our surface, we follow the main lines of H. Karcher's
conjugate Plateau construction. The principle is to first construct
a minimal surface with boundary made of straight edges. Typically
this surface is constructed by solving a Dirichlet boundary value
problem on a polygonal domain with piecewise linear boundary data.
(The boundary data may take on infinite values, in which case this is
called a Jenkins Serrin type problem.)
Then one considers the conjugate minimal
surface which is bounded by planar geodesics.
If the polygonal boundary of the first surface is well chosen, the
conjugate surface will extend by symmetry to an embedded, complete minimal
surface.
Adjusting the lengths of the edges so that this is the case is called the
Period Problem. The difficulty of solving this problem is the main limitation
of the method.

\medskip

We add one more step to this construction. We first solve a Dirichlet boundary
value problem for the {\em maximal} graph equation, with piecewise linear
boundary data. Then we consider the conjugate function, whose graph is
minimal. Then as above we consider the conjugate minimal surface, solve the
Period Problem and extend by symmetry. 
So there are two consecutive conjugations, altough of a different nature.
The advantage of this approach is that part of the Period Problem (namely the
vertical part) will be automatically solved.
More details on maximal surfaces will be given in section
\ref{section-preliminaires}. 

\medskip

In our case, since we add infinitely many handles in a non-periodic
way, we are faced with an infinite dimensional Period Problem.
We begin by adding a finite number $N$ of handles.
We solve the Period
Problem using the Poincar\'e Miranda Theorem, which is a natural 
$N$-dimensional extension of the intermediate value theorem.
Then we let $N\to\infty$.

\section{Preliminaries}
\label{section-preliminaires}
\subsection{Minimal and maximal graphs}
Let $u$ be a function on a domain $\Omega\subset\R^2$. The graph of
$u$ is a minimal surface if $u$ satisfies the minimal graph equation
\begin{equation}
\label{eq-min}
\mbox{div}\left(\frac{\nabla u}{\sqrt{1+|\nabla u|^2}}\right)=0
\end{equation}
This equation is equivalent to the fact that the conjugate 1-form
$$d\Psi_u=\frac{u_x}{\sqrt{1+|\nabla u|^2}}dy 
-\frac{u_y}{\sqrt{1+|\nabla u|^2}}dx$$
is closed.
Locally, $d\Psi_u$ is then the differential of a function $\Psi_u$
called the conjugate function.
Then $v=\Psi_u$ is a solution of the maximal graph equation
\begin{equation}
\label{eq-max}
\mbox{div}\left(\frac{\nabla v}{\sqrt{1-|\nabla v|^2}}\right)=0
\end{equation}
This is called the maximal graph equation because $v$ satisfies
(\ref{eq-max}) if its graph is a maximal surface in the 
Lorentzian space $\L^3$, namely a space-like surface which is
a critical point for the area functionnal.

\medskip

Conversely, let $v$ be a solution of (\ref{eq-max}). Then the
conjugate 1-form
$$d\Phi_v=\frac{v_y}{\sqrt{1-|\nabla v|^2}}dx 
-\frac{v_x}{\sqrt{1-|\nabla v|^2}}dy$$
is closed. Hence locally $d\Phi_v$ is the differential of a function
$u=\Phi_v$ which solves the minimal graph equation (\ref{eq-min}).
Moreover up to a constant, $\Phi_{\Psi_u}=u$.
\subsection{The Dirichlet boundary value problem}
Let $\Omega\subset\R^2$ be a bounded domain.
Let $v:\Omega\to\R$ be a smooth function satisfying (\ref{eq-max}).
Then $|\nabla v|<1$ hence $v$ is Lipschitz and extends continuously
to $\partial\Omega$, so
we can talk about the boundary values of $v$.
(For this to be true, we need some regularity of the boundary of
$\Omega$. All the domains we consider will have piecewise smooth
boundary.)

\medskip
We need to construct solutions $v$ of the maximal graph equation
(\ref{eq-max}) in $\Omega$, with precribed boundary values, and
with singularities at some prescribed points inside $\Omega$.
For this we use the following theorem, which is a consequence
of Theorem 1 in \cite{klyachin-miklyukov}
and Theorem 4.1 in \cite{bartnik-simon} :
\begin{theorem}
\label{th-max}
Let $\Omega\subset\R^2$ be a bounded domain.
Let $\X\subset\Omega$ be a finite set (the singular set).
Let $\varphi : \partial\Omega\cup \X\to \R$ be a given function such that
\begin{equation}
\label{eq-condition-max}
\forall p,p'\in\partial\Omega\cup \X,\quad p\neq p',\quad
|\varphi(p)-\varphi(p')|\leq d_{\Omega}(p,p')
\end{equation}
where the inequality is strict whenever the segment $[p,p']$ is not
contained in $\partial\Omega$.
Then there exists a function $v:\Omega\to\R$
which satisfies the maximal graph equation (\ref{eq-max})
in $\Omega\setminus \X$, with boundary data
$v=\varphi$ on $\partial\Omega\cup \X$.
This function is smooth in $\Omega\setminus \X$.
(Here $d_{\Omega}$ is the intrinsic distance of $\Omega$, so if
$\Omega$ is convex, it agrees with the euclidean distance.)
\end{theorem}

\subsection{Some complements on the correspondence $v\leftrightarrow
  \Phi_v$}

Since we obtain solutions $v$ to the Dirichlet problem for
\eqref{eq-max}, we need to understand the behaviour of the conjugate
function $\Phi_{v}$ near the boundary.
The first result describes the behaviour near the boundary of the
domain $\Ome$.
\begin{lemma}[\cite{jes1} and \cite{mazet2}]
Let $v$ be a solution of \eqref{eq-max} on $\Ome$ and $T\subset\partial
\Ome$ be an open straight segment oriented as $\partial \Ome$. Then $\int_T
dv=|T|$ if and only if $\Phi_v$ diverges to $+\infty$ on $T$.
\end{lemma}
Now we shall describe the behaviour near a singularity in the
domain. Let $v$ be a solution of \eqref{eq-max} on a punctured disk
$\D^*$ (with $\D=\{(x,y)\in\R^2|\,x^2+y^2<r^2\}$). Then the conjugate
function $u=\Phi_v$ is not well defined on 
$\D^*$; actually, $u$ is multivalued in the sense that when we turn
around the origin we need to add a constant to $u$: this constant is
given by $\int_\gamma d\Phi_v$ where $\gamma$ generates
$\pi_1(\D^*)$. If this constant vanishes, $u$ is well defined and so extends
smoothly to the whole disk; $v$ then also extends to $\D$ and the origin is
a removable singularity for $v$.

In the case $\int_\gamma d\Phi_v\neq0$, the graph of the multivalued
function $u$ has then the shape of a half-helicoid. On the boundary of
the cylinder $\D\times\R$ the graph is bounded by a helix-like looking
curve. It is bounded by a vertical straight line above the origin.

In fact in the paper, we are always in the case where $v$ is positive
and vanishes at the origin. This first implies that $\int_\gamma
d\Phi_v\neq0$. Besides, Theorem 4.2 in \cite{mazet2} proves that the
graph of $u$ is bounded by a vertical straight line above the origin.

\subsection{Convergence of sequences of solutions}
We shall study many times the convergence or the divergence of
sequences of solutions to \eqref{eq-max}. In this subsection, we
expose some results that we will use. Actually, these
results were developped by the first author in \cite{mazet1,mazet2} for
solutions
of \eqref{eq-min}; the correspondence $u\leftrightarrow \Psi_u$ and
$v\leftrightarrow \Phi_v$ translates them to solutions of
\eqref{eq-max}. Here the convergence that we shall consider is the
$C^k$ convergence on compact subsets of the domain for every $k$.

So let us consider a sequence $(v_n)_{n\in\N}$ of solutions to
\eqref{eq-max} which are defined on a domain $\Ome$. We first notice that
since each $v_n$ is Lipschitz continuous there exists a subsequence of
$(v_n-v_n(q))$ (where $q\in \Ome$ is a fixed point)
that converges to a lipschitz function $v$
on $\Ome$; but the convergence is only $C^0$ and $v$ is \textit{a priori}
not a solution to \eqref{eq-max}. However, since we have convergence on
$\overline{\Ome}$, we can talk about the boundary value of the limit.

To study the smooth convergence, we first define the \emph{convergence
  domain} of the sequence by 
$$\boB((v_n)_{n\in\N})=\{p\in \Ome\,|\, \sup_n \{|\nabla v_n|(p)\}<1\}$$

$\boB((v_n)_{n\in\N})$ is an open subset of $\Ome$ and on each component
$\Ome'$ of it, there is a subsequence of $(v_n-v_n(q))_{n\in\N}$
converging $C^k$ on compact subsets of $\Ome'$ to a solution $v$ of
\eqref{eq-max}, where $q$ is some fixed point in $\Ome'$. We notice
that all solutions of \eqref{eq-max} that we shall consider are bounded by
$1$; thus we do not need to use the vertical translation by $v_n(q)$ to
ensure the convergence. Besides $\Ome\setminus\boB((v_n)_{n\in\N})$ is
the union of straigt lines 
$\cup_iL_i$, where each $L_i$ is a component of the intersection of a
straight line with $\Ome$. The $L_i$ are called \emph{divergence lines} of
the sequence $(v_n)_{n\in\N}$ since $\sup_n \{|\nabla v_n|(p)\}=1$ for
$p\in L_i$; more precisely we have
\begin{lemma}
Let $p$ be a point in a divergence line $L$, then a subsequence of
$(\nabla v_n(p))_{n\in\N}$ converges to one of the two unit generating
vectors of $L$.

Besides, if $T$ is a segment in $L$, it holds $\int_T dv_n\rightarrow
\pm|T|$ for a subsequence. 
\end{lemma}

To ensure the convergence of a subsequence of $(v_n)_{n\in\N}$ on $\Ome$
it then suffices to prove there are no divergence line. The above lemma
is one tool in that direction. The following one is another.
\begin{lemma}\label{lem-div-bdry}
Let us assume that one part of the boundary of $\Ome$ is a segment
$[a,b]$. If for every $n$ $|v_n(a)-v_n(b)|=|ab|$, then no divergence
line can end in the interior of $[a,b]$.
\end{lemma}

Actually in this paper, the solutions are not defined on the same domain
$\Ome$: we have in fact a sequence of domains $(\Ome_n)_{n\in\N}$ and each
solution $v_n$ is defined on $\Ome_n$. So to make sense to the above
definition we need to introduce the \emph{limit domain} $\Ome_\infty$:
$$\Ome_\infty=\bigcup_{p\in\N} \text{Int}\left(\bigcap_{k\ge p}\Ome_k
\right)$$ 
A point is then in $\Ome_\infty$ if a neighborhood of this point is
included in all $\Ome_k$ for $k$ great enough. With this definition, we
have anew the convergence domain and the divergence lines by replacing
$\Ome$ by $\Ome_\infty$.

We notice that when $(\Ome_n)_{n\in\N}$ is an increasing sequence,
$\Ome_\infty$ is simply the union of all the $\Ome_n$. In this paper,
the sequence $\Ome_n$ is often $\Ome\setminus\boS_n$ where $\boS_n$ is
a locally finite set of points. If $(\boS_n)$ converges on compact 
subsets to a locally finite subset $\boS_\infty$ then
$\Ome_\infty=D\setminus \boS_\infty$.

\section{The fundamental piece}
\subsection{The Dirichlet boundary value problem}
In this section we solve a Dirichlet boundary value problem for the
maximal graph equation (\ref{eq-max}) in an infinite strip.
The solution $v$ will have singularities at some prescribed points.
The position of these singularities are the parameters of our construction.
(Each singularity is responsible for one handle of the minimal surface we are
constructing.
In later sections, we will adjust these parameters so that the Period
Problem is solved.)

\medskip

\begin{figure}[h]
\begin{center}
\resizebox{0.8\linewidth}{!}{\input{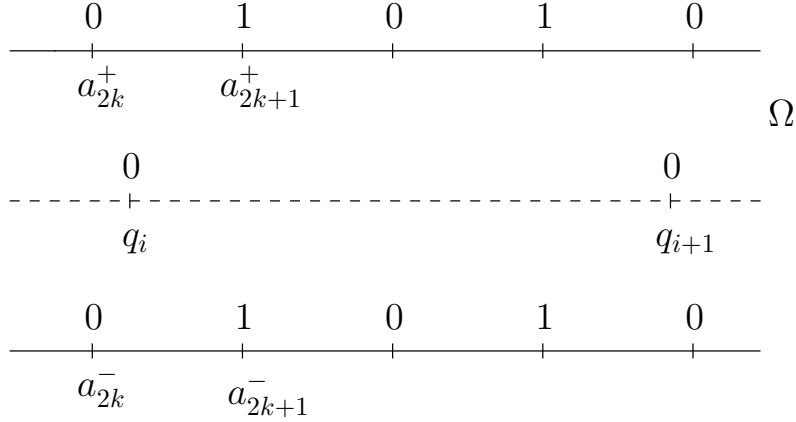}}
\caption{the Dirichlet boundary value problem.}
\label{fig1}
\end{center}
\end{figure}

Fix some $\ell>0$ and let $\Omega$ be the strip 
$\R\times(-\ell,\ell)$.
Let us define the boundary data $\varphi$ on $\partial\Omega$ as
follows :
for $k\in \Z$, let $a^+_k=(k,\ell)$ and $a^-_k=(k,-\ell)$.
Define $\varphi$ on the segment $[a^{\pm}_{2k-1},a^{\pm}_{2k+1}]$ by
$\varphi(p)=|p-a^{\pm}_{2k}|$.
In other words $\varphi$ is piecewise affine on $\partial\Omega$,
with value $0$ at $a^{\pm}_{2k}$ and $1$ at $a^{\pm}_{2k+1}$
(see figure \ref{fig1}).

\medskip

Let $\X$ be a closed, discrete subset of the horizontal line $y=0$. It
will be convenient to identify the $x$-axis with $\R$ and see $\X$ as
a subset of $\R$. When $\X$ is finite, we write
$\X=\{q_1,\cdots,q_N\}$ and assume that $q_1<q_2<\cdots<q_N$.
When $\X$ is infinite, we may write $\X=\{q_i : i\in I\}$, with
$q_i<q_{i+1}$, where $I$ is either $\N$, $-\N$ or $\Z$, depending
on whether $\X$ is bounded from below, bounded from above, or neither.
Finally, we define $\varphi=0$ on $\X$.

\begin{proposition}
\label{prop-max}
Let $\Omega$ and $\X$ be as above.
Assume that 
\begin{equation}
\label{eq-condition}
\forall q\in \X,\quad\forall k\in\Z,\quad
|q-a^+_{2k+1}|>1.
\end{equation}
Then there exists a function $v$ on $\Omega$ which
solves (\ref{eq-max}) in $\Omega\setminus \X$,
with boundary data $v=\varphi$ on $\partial
\Omega\cup \X$.
Moreover, $0\leq v\leq 1$ in $\Omega$.
The function $v$ is unique.
\end{proposition}
We call the function $v$ the solution to the
Dirichlet problem in $\Omega\setminus\X$.
When needed, we will write $v=v[\X]$.
\begin{remark}
The condition (\ref{eq-condition}) is automatically
satisfied when $\ell>1$. We are however mostly interested in the
case $\ell<1$, as this is the only case where we know how to solve
the Period Problem.
\end{remark}
Proof of the proposition : for $n\in\N^*$,
consider the box $\Omega_n=(-2n,2n)\times(-\ell,\ell)$.
We first solve the Dirichlet problem on $\Omega_n$
and then let $n\to\infty$.
Let $\X_n=\X\cap\Omega_n$.
Define $\varphi_n$ on $\partial\Omega_n\cup \X_n$ by
$\varphi_n=\varphi$ on the horizontal edges
$[-2n,2n]\times\{\pm\ell\}$,
$\varphi_n=0$ on the vertical edges $\{\pm 2n\}\times[-\ell,\ell]$,
and $\varphi_n=0$ on the singular set $\X_n$.
\begin{claim}
\label{claim-max}
The function $\varphi_n$ on $\partial \Omega_n\cup \X_n$
satisfies the condition (\ref{eq-condition-max}) of theorem \ref{th-max}.
\end{claim}
Proof :
consider $p,p'\in\partial\Omega_n\cup \X_n$, $p\neq p'$.
\begin{itemize}
\item If $p$ and $p'$ are both on the line $y=\ell$, then clearly
$|\varphi_n(p)-\varphi_n(p')|\leq |p-p'|$.
\item
If $p$ and $p'$ are on $\partial\Omega_n$,
let $\widetilde{p}$ and $\widetilde{p}'$ be the projections of $p$ and $p'$
on the line $y=\ell$. Then 
$$|\varphi_n(p)-\varphi_n(p')|= |\varphi_n(\widetilde{p})-
\varphi_n(\widetilde{p}')| 
\leq |\widetilde{p}-\widetilde{p}'|\leq |p-p'|.$$
Moreover, if the segment $[p,p']$ is not horizontal, the
last inequality is strict. If the segment $[p,p']$ is
horizontal, and is not included in $\partial\Omega_n$,
then $p$ and $p'$ are both on the vertical edges,
so $|\varphi_n(p)-\varphi_n(p')|=0<|p-p'|$ as required.
\item If $p$ is on $\partial\Omega_n$ and $p'=q\in \X$ :
if $p$ is on a vertical edge, then $\varphi_n(p)=\varphi_n(q)=0$.
If $p$ is on the segment $[a^+_{2k},a^+_{2k+1}]$, we have
$$|\varphi_n(p)-\varphi_n(q)|=|p-a^+_{2k}|
=1-|p-a^+_{2k+1}|\leq 1+|q-p|-|q-a^+_{2k+1}|<|p-q|$$
where we have used the triangle inequality and the hypothesis of
proposition \ref{prop-max}.
The case where $p$ is on the segment $[a^+_{2k-1},a^+_{2k}]$ is similar,
and the case where $p$ is on the line $y=-\ell$ follows by symmetry of
$\varphi_n$.
\item If $p$, $p'$ are both in $\X$, then $\varphi_n(p)=\varphi_n(p')=0$.
\end{itemize}
\cqfd

By theorem \ref{th-max}, there exists a solution $v_n$ of the
maximal graph equation (\ref{eq-max}) on $\Omega_n\setminus \X_n$
with boundary data $\varphi_n$.
Since $v_n$ extends continously to the compact set
$\overline{\Omega_n}$, $v_n$ is bounded.
By the maximum principle for the maximal graph equation,
$v_n$ reaches its maximum and its minimum at a boundary point or
a singular point, so $0\leq v_n\leq 1$ in $\Omega_n$.
Consider now the sequence $(v_n)_n$. Let $L$ be a divergence line.
Let $T\subset L$ be
a segment, then $\lim\int_T dv_n =\pm |T|$. Since $v_n$ is bounded, this
implies that $L$ has finite length so $L$ is a segment connecting two
points $p$ and $p'$ on $\partial\Omega\cup \X$. Then
$$|\varphi(p)-\varphi(p')|=
|\int_p^{p'} dv_n|\to|p-p'|\quad\Rightarrow\quad
|\varphi(p)-\varphi(p')|=|p-p'|$$
which contradicts claim \ref{claim-max} since $L\subset\Omega$.
Hence there are no divergence lines, so passing to a subsequence,
$(v_n)_n$ converges on compact subsets of $\overline{\Omega}$
to a function $v$, which is a solution of (\ref{eq-max}) in
$\Omega\setminus \X$ with boundary data $\varphi$ on
$\partial\Omega\cup \X$. Uniqueness follows from Theorem 2 in
\cite{mazet3}. 
\cqfd

\subsection{The minimal graph}
In this subsection and the following one, we assume that
$\X=\{q_1,\cdots,q_N\}$ with $q_1<\cdots<q_N$.
Let $v$ be the solution of the Dirichlet problem on
$\Omega\setminus\X$, given by proposition \ref{prop-max}.
Let $\Omega^+$ be the strip $\R\times(0,\ell)$.
Since $\Omega^+$ is simply connected and $v$ is smooth in $\Omega^+$,
the conjugate function $u$ is well defined (up to a constant) in
$\Omega^+$.
The graph of $u$ is a minimal surface.
In this section we describe geometrically its boundary.

\medskip

By uniqueness, $v$ satisfies $v(x,-y)=v(x,y)$ in $\Omega$.
Hence on the $x$-axis,
away from the singular points $q_1,\cdots,q_N$,
we have $v_y=0$.
From the definition of $\Phi_v$, this gives
$u_x=0$. Hence $u$ is locally constant on the $x$-axis minus
the singular points,
with a finite number of jumps at the points $q_1,\cdots,q_N$.
On the line $y=\ell$, $u$ goes to $+\infty$ on the segments 
$(a^+_{2k-1},a^+_{2k})$ and to $-\infty$ on the
segments $(a^+_{2k},a^+_{2k+1})$, $k\in\Z$.

\medskip

\begin{figure}[h]
\begin{center}
\resizebox{0.8\linewidth}{!}{\input{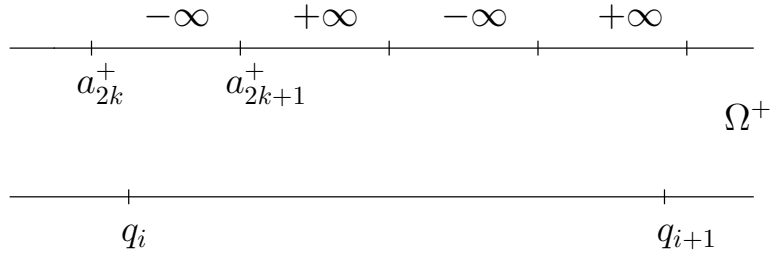}}
\caption{the minimal graph}
\label{fig2}
\end{center}
\end{figure}

Let $M$ be the graph of $u$ on the strip $\Omega^+$. The minimal surface $M$ 
is bounded by vertical lines $A_k$ above the points $a^+_k$, $k\in\Z$,
by $N$ vertical segments $B_i$ above the points $q_i$, $i=1,\cdots,N$,
by $N-1$ horizontal segments $C_i$ above the segments $(q_i,q_{i+1})$,
$i=1,\cdots,N-1$
and by two horizontal half-lines $C_0$ and $C_N$ above
$(-\infty,q_1)$ and $(q_N,+\infty)$ (see figure \ref{fig2}).
The heights of the horizontal pieces $C_0,\cdots,C_N$ are unknown.

\subsection{The conjugate minimal surface}
\begin{figure}[h]
\begin{center}
\resizebox{0.7\linewidth}{!}{\input{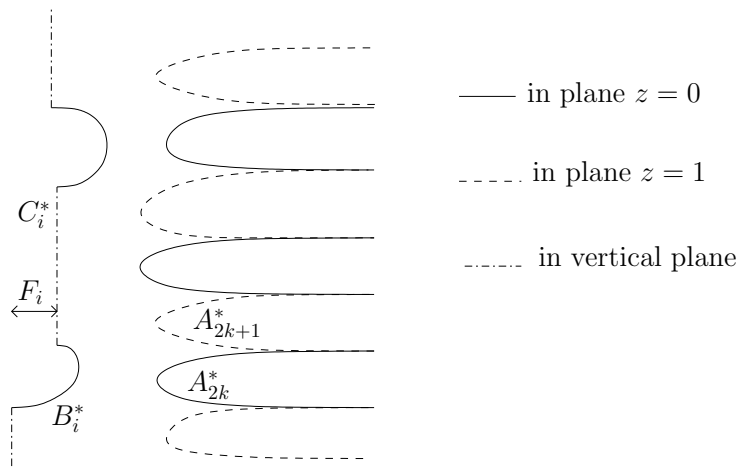}}
\caption{the conjugate minimal surface, Period Problem not solved}
\label{fig3}
\end{center}
\end{figure}

Let $M^*$ be the conjugate minimal surface to $M$. The third
coordinate of $M^*$ (seen as an immersion of the strip $\Omega^+$) is
the function $v$, so $M^*$ lies in the slab $0\leq z\leq 1$. Let
$A_k^*$, $B_i^*$ and $C_i^*$ denote the corresponding conjugate curves
on $M^*$ (see figure \ref{fig3}).
Then the $A_k^*$, $k\in\Z$, and $B_i^*$, $i=1,\cdots,N$, are
horizontal geodesics. From the boundary values of $v$, $A_{2k}^*$ and
$B_i^*$ lie in the plane $z=0$, while $A_{2k+1}^*$ lies in the plane
$z=1$. Each $C_i^*$, $i=0,\cdots,N$ is a geodesic contained in a
vertical plane parallel to the plane $x=0$. There is no reason however
that all $C_i^*$ are in the same vertical plane : this is the Period
Problem, which we will consider in the next section.

\begin{remark}
The function $u$ is the solution of a Jenkins Serrin type problem on
the strip $\Omega^+$.
One possible way to construct $M$ would be to directly solve
this Jenkins Serrin problem, with the jumps of $u$ at each
point $q_i$ as parameters.
Then we would have to adjust these parameters to guarantee that
the conjugate curves $B_i^*$ all lie in the plane $z=0$.
This means another Period Problem to solve. It is automatically
solved in our maximal graph approach.
\end{remark}
\section{The Period Problem}
In this section, we first formulate the Period Problem in
general. Then we solve it provided that $\ell<1$ and there is a finite
number of singularities $q_1,\cdots,q_N$, which are not too close
from each other. Our solution to the Period Problem is based on the
Poincar\'e Miranda theorem :
\begin{theorem}[Poincar\'e Miranda]
Let $F=(F_1,\cdots,F_N)$ be a continous map from $[0,1]^N$
to $\R^N$.
Write $x=(x_1,\cdots,x_N)$.
Assume that for each $i$, $F_i(x)$ is negative on
the face $x_i=0$, while $F_i(x)$ is positive on the face
$x_i=1$. Then there exists $x^0\in [0,1]^N$ such that
$F(x^0)=0$.
\end{theorem}

\subsection{Formulation of the Period Problem}
Let $\Omega$ and $\X$ be as in proposition \ref{prop-max}, and
let $v$ be the solution of the Dirichlet problem in $\Omega\setminus\X$.
Let $u$ be the conjugate function of $v$ and $X^*=(X_1^*,X_2^*,X_3^*)$
be the conjugate minimal surface to the graph of $u$. Both $u$ and $X^*$
are only locally well defined, but their differentials are well defined
in $\Omega\setminus \X$.
Explicitely, $dX^*$ is given by
\begin{equation}
\label{eq-dX1star}
dX_1^*=\frac{u_x u_y dx + (1+(u_y)^2)dy}{\sqrt{1+|\nabla u|^2}}
\end{equation}
\begin{equation}
\label{eq-dX2star}
dX_2^*=\frac{ - (1+(u_x)^2)dx - u_x u_y dy}{\sqrt{1+|\nabla u|^2}}
\end{equation}
$$dX_3^*=dv$$
The Period problem asks that $X^*$ is well defined in $\Omega\setminus \X$.
This is equivalent to $\int_{\gamma}dX^*=0$, where $\gamma$ is a small
circle around any point of the singular set $\X$.
\begin{proposition}
\label{prop-symmetries}
$X_2^*$ and $X_3^*$ are well defined in $\Omega\setminus\X$.
\end{proposition}
Proof : this is clear for $X_3^*$. For $X_2^*$, we use the following symmetry
argument.
Let $\tau(x,y)=(x,-y)$.
By uniqueness and symmetry of the boundary data, $v\circ\tau=v$.
Hence $\tau^* dv =dv$, so $v_x\circ\tau=v_x$ and
$v_y\circ\tau=-v_y$.
This gives $u_x\circ\tau=-u_x$ and $u_y\circ\tau=u_y$.
Hence using (\ref{eq-dX2star}), $\tau^* dX_2^*=dX_2^*$.
Let $\gamma$ be a small circle around a singularity $q\in\X$.
Since $\tau(\gamma)$ is homologous to $-\gamma$, this
gives $\int_{\gamma}dX_2^*=0$, so $X_2^*$ is well defined in 
$\Omega\setminus\X$.
This also gives $\tau^*dX_1^*=-dX_1^*$, so $X_1^*$ is locally constant on
the $x$-axis.
\cqfd

\medskip

By proposition  \ref{prop-symmetries},
we only have to worry about the periods of $dX_1^*$.

From now on, we assume that $\X=\{q_1,\cdots,q_N\}$ is finite.
Let $\gamma_i$ be a small circle around the point $q_i$ and let
$$F_i(q_1,\cdots,q_N)=\int_{\gamma_i} dX_1^*.$$
The Period Problem asks that $F_i=0$ for $1\leq i \leq N$.
Note that by symmetry, $F_i$ is equal to twice the integral of
$dX_1^*$ on a half circle from $q_i+\epsilon$ to $q_i-\epsilon$,
so $F_i=0$ means that the curves $C_i^*$ and $C_{i-1}^*$
are in the same vertical plane as required (see figure \ref{fig3}).
\subsection{Continuity of the periods}
To apply the Poincar\'e Miranda theorem, 
we need the continuity of the periods with respect to the parameters.
\begin{proposition}
\label{prop-continuity}
The periods $F_i$ depend continuously on $(q_1,\cdots,q_N)$.
\end{proposition}
Proof : consider an admissible value $(q_1,\cdots,q_N)$ of the
parameters (namely, all $q_j$ satisfy equation (\ref{eq-condition})).
Consider a sequence $(q^n_1,\cdots,q^n_N)$ converging to
$(q_1,\cdots,q_N)$. Let $\X_n=\{q^n_1,\cdots,q^n_N\}$ and
$\X=\{q_1,\cdots,q_N\}$. Let $v_n$ and $v$ be the solutions of the
Dirichlet problem in $\Omega\setminus\X_n$ and $\Omega\setminus\X$,
respectively. Assume the sequence $(v_n)_n$ has a divergence line.
Then arguing as in the proof of proposition \ref{prop-max},
$L$ has finite length so is a segment connecting two points
of $\partial\Omega\setminus\X$, which contradicts in the same way the
fact that the points $q_j$ satisfy (\ref{eq-condition}).
Hence there are no divergence lines, so a subsequence of $(v_n)_n$
converges on compact subsets of $\overline{\Omega}$ to a solution to the
Dirichlet problem on $\Omega\setminus\X$. By uniqueness of the
solution to this problem, the whole sequence $(v_n)_n$ converges to
$v$ on compact subsets of $\overline{\Omega}$, and converges smoothly
to $v$ on compact subsets of $\Omega\setminus\X$. This implies that
$du_n$ converges to $du$ and $dX_{1,n}^*$ converges to $dX_1^*$ on
compact subsets of $\Omega\setminus\X$. Integrating on $\gamma_i$
which encloses $q_i^n$ for big $n$, we obtain that
$F_i(q_1^n,\cdots,q_N^n)\to F_i(q_1,\cdots,q_N)$.  
\cqfd
\subsection{Local property of the period}
From now on we assume that $\ell<1$.
Let us define for the rest of the paper
$$\eta=1-\sqrt{1-\ell^2}.$$
In this section, we prove that some properties of the period
$F_i(q_1,\cdots, q_N)$ depends only on the position of $q_i$ if the
other $q_j$ are not too close from $q_i$.

Let us denote by $\Ome_L$ the box $(-L,L)\times(-\ell,\ell)$ and
consider $q\in(-\eta,\eta)$ and a finite set of points $\boS$ in
$(-2,-2+\eta)\cup (2-\eta,2)$. Let $v$ be a solution of the maximal
graph equation \eqref{eq-max} on $\Ome_2\backslash(\{q\}\cup\boS)$, with
boundary value $v=\phi$ on $(-2,2)\times\{-\ell,\ell\}$, $v(q)=0$ and
$v(\boS)=0$. The boundary value of $v$ on the vertical edges is free,
although we require $0\le v\le 1$. Let us study the divergence lines
of a sequence of such solutions $v$.
\begin{lemma}\label{convergence}
For every $n\in\N$, let $q_n$, $\boS_n$ and $v_n$ be as above. We assume that
$\lim q_n=q$ exists. Then 
\begin{itemize}
\item if $q\in(-\eta,\eta)$, there is no divergence line in
  $\Ome_1\setminus\{q\}$. 
\item if $q=\eta$, the only divergence lines in
  $\Ome_1\setminus\{\eta\}$ are $[\eta,a_1^-]$ and $[\eta,a_1^+]$.
\end{itemize}
\end{lemma}
Proof: 
The two segments $[\eta,a_1^+]$ and $[\eta,a_1^-]$ have length
one. Since $v_n(q_n)=0$ and $v_n(a_1^+)=1=v_n(a_1^-)$, both segments
are divergence lines for $(v_n)_{n\in\N}$ in the second case.

Let us now prove that there is no other divergence line in $\Ome_1$
for both cases. Since $0\le v_n\le 1$ every divergence line is a
segment of length at most one. Hence a divergence line $L$ which
intersects $\Ome_1$ must have an end-point in $\overline{\Ome_1}$.
Because of 
lemma \ref{lem-div-bdry} these end-points needs to be $a_0^\pm$,
$a_{-1}^\pm$, $a_1^\pm$ or $q$. Let us assume that $a_0^+$ is one
end-point of $L$. The distance from $a_0^+$ to $\Ome_2\setminus\Ome_1$
is one, hence the other end-point is in $\Ome_1$. It can not be $a_0^-$
or $q$ since $v_n(a_0^+)=v_n(a_0^-)=\lim v_n(q)$. It is not
$a_{\pm1}^-$ since the distance from $a_0^+$ to these points is
$\sqrt{1+4\ell^2}>1$. Then $a_0^+$ is not an end point of $L$; by
symmetry, this is also true for $a_0^-$. Let us assume that $a_{-1}^+$
is an end-point of $L$ then the other end-point is either $q$ or
$a_1^-$ but the distance from $a_{-1}^+$ to these two points is
strictly larger than one since $q\in(-\eta,\eta]$; then $a_{-1}^+$
  is not an end point for $L$. By symmetry, this is also true for
  $a_{-1}^-$ and $a_1^\pm$ unless $q=\eta$ and
  $L=[\eta,a_1^\pm]$ which is the case we studied first. Then we
  can assume that $q$ is an end-point of $L$ and the other one is in
  $\Ome_2\setminus\Ome_1$. By lemma \ref{lem-div-bdry}, if $L$ is not
  horizontal, the second end-point is on the vertical edges of
  $\Ome_2$ but the distance from $q$ to these edges is larger than
  $2-\eta>1$. So $L$ is horizontal; we assume, for example, that $L$
  is on the left of $q$. Since the length of $L$ is less than $1$, the
  other end-point of $L$ needs to be in the interior of $\Ome_2$ and
  $\boS_n\cap (2-\eta,2]\neq\emptyset$. Let $s_n =\min \boS_n\cap
    (2-\eta,2]$. We assume that $(s_n)$ converges to $s$ in
      $[2-\eta,2]$, then $L$ is the segment $[q,s]$. We have
      $v_n(s_n)=0$. Since $L$ is a divergence line $|L|=\lim
      |v_n(q)-v_n(s)|=\lim |v_n(q_n)-v_n(s_n)|=0$ this gives a
      contradiction and the lemma is proved. \cqfd 

Let $v$ be a solution of the maximal
graph equation \eqref{eq-max} on $\Ome_2\backslash(\{q\}\cup\boS)$ as
above, besides we assume now the symmetry $v(x,y)=v(x,-y)$. In
applications, $v$ will be the restriction of some $v[q_1,\dots,q_n]$
to a box around one $q_i$. Let $\gamma$ be a small circle around $q$
then we define the period $F(v)$ by $\int_\gamma d X_1^*$ where $d
X_1^*$ is given by equation \eqref{eq-dX1star} with $u$ the conjugate
function to $v$. We then have some control on the behaviour of the period.
\begin{proposition}\label{local-bhv-period}
There exists $\eta_0\in(0,\eta)$ which depends only on $\ell$ such
that for any solution $v$ of the above Dirichlet problem on
$\Ome_2\backslash(\{q\}\cup \boS)$ we have:
\begin{itemize}
\item if $\eta_0\le q<\eta$, $F(v)\ge 1$.
\item if $-\eta<q\le -\eta_0$, $F(v)\le-1$.
\end{itemize}
\end{proposition}
Proof : Let $\sigma(x,y)=-(x,y)$ and $v$ be a solution of the above
Dirichlet problem on $\Ome_2\backslash(\{q\}\cup\boS)$ then $v'=v\circ
\sigma$ is a solution of this Dirichlet problem on
$\Ome_2\backslash(\{-q\}\cup -\boS)$. From the definition of
$du=d\Phi_v$ and $du'=d\Phi_{v'}$, we obtain $\sigma^* du =du'$. From
equation (\ref{eq-dX1star}), we get $\sigma^* dX_1^*=-{dX_1^*}'$ where
$dX_1^*$ and ${dX_1^*}'$ are respectively associated to $v$ and
$v'$. Since $\sigma$ preserves orientation, integrating on a small
circle around $q$ gives $F(v)=-F(v')$. Thus the second item of the
proposition is a consequence of the first one. 

If the first item is wrong there exists a sequence $q_n\rightarrow
\eta$ and for each $n$ a set $\boS_n$ and a solution $v_n$ of the
above Dirichlet problem on $\Ome_2\setminus(\{q_n\}\cup\boS_n)$ such
that $F(v_n)<1$. Let us prove that, actually, $\lim F(v_n)=\infty$.

By lemma \ref{convergence} the two segments $L^+=[\eta,a_1^+]$ and
$L^-=[\eta,a_1^-]$ are divergence lines for $(v_n)_{n\in\N}$. On
$L^+$, $(\nabla v_n)_{n\in\N}$ converges to $\overrightarrow{\eta
  a_1^+}=(\sqrt{1-\ell^2},\ell)$. Let $\Ome_-$ be the connected
component of $\Ome_1\backslash(L^+\cup L^-)$ containing the
origin. Because of lemma \ref{convergence}, $\Ome_-$ is included in
the convergence domain of $(v_n)_{n\in\N}$ 

Then we can assume that the sequence $(v_n)_{n\in\N}$ converges on
$\Ome_-$ to a solution $v$ which takes on the boundary the value
$\phi$ on $\partial\Ome \cap \Ome_-$ and $|p-\eta|$ for $p\in
L^+\cup L^-$. For every $n$, let $u_n$ be the conjugate function
$\Phi_{v_n}$ which is defined on $\Ome_-\backslash[q_n,\eta]$. The
limit domain of $(\Ome_-\backslash[q_n,\eta])_{n\in\N}$ is $\Ome_-$
and $(u_n)_{n\in\N}$ converges to $u=\Phi_{v}$. Because of the
boundary value of $v$, $u$ takes the value $+\infty$ along $L^+$. We
are interested in what happens near the middle point
$((1+\eta)/2,\ell/2)$ of $L^+$. Because of Lemma 1 in \cite{jes1} we
have : 
\begin{equation}\label{lim-u_y}
\frac{{u}_y}{\sqrt{1+|\nabla u|^2}}(\frac{1+\eta}{2}, y)
  \longrightarrow -\sqrt{1-\ell^2}
\end{equation}
when $y\rightarrow \ell/2$ with $y>\ell/2$. Lemma 1 in \cite{jes1}
  implies also that:
\begin{equation}\label{equiv-u_y}
{u}_y(\frac{1+\eta}{2}, y)\le -\frac{C}{|y-\ell/2|},\quad (C>0)
\end{equation}
for $y>\ell/2$ near $\ell/2$.

Let us consider the following path $\Gamma$: it is the union of the
segment $[((1+\eta)/2,0),((1+\eta)/2,3\ell/4)]$ with a curve
$\Gamma_3$ in $\Ome_-\cap\{y>0\}$ that joins $((1+\eta)/2,3\ell/4)$ to 
$(0,0)$. For $t\in(\ell/2,3\ell/4)$, let $\Gamma_1(t)$ be the segment
$[((1+\eta)/2,0),((1+\eta)/2,t)]$ and $\Gamma_2(t)$ be the segment
$[((1+\eta)/2,t),((1+\eta)/2,3\ell/4)]$. 

\begin{figure}[h]
\begin{center}
\resizebox{0.6\linewidth}{!}{\input{figquasi1.pstex_t}}
\caption{}
\end{center}
\end{figure}

Because of the symmetry $v_n(x,y)=v_n(x-y)$, for large $n$, the period
$F(v_n)$ is given by $2\int_{\Gamma}{dX_1^*}_{,n}$ where
${dX_1^*}_{,n}$ is associated to $u_n$. Because of
Equation \eqref{eq-dX1star}, $\int_{\Gamma_1(t)}{dX_1^*}_{,n}\ge 0$. Hence
\begin{equation}\label{minor-period}
F(v_n)\ge 2\int_{\Gamma_2(t)\cup\Gamma_3}{dX_1^*}_{,n} 
\end{equation} 
Since the convergence $u_n\rightarrow u$ is smooth on compact
subsets of $\Ome_-$, 
\begin{equation}\label{limit-G2G3}
\int_{\Gamma_2(t)\cup\Gamma_3}{dX_1^*}_{,n} \longrightarrow
\int_{\Gamma_2(t)\cup\Gamma_3}dX_{1}^*
\end{equation}
with $dX_{1}^*$ associated to $u$. By \eqref{lim-u_y} and
\eqref{equiv-u_y}, we have 
\begin{equation}\label{limit-G2t}
\int_{\Gamma_2(t)}dX_{1}^*\xrightarrow[t\rightarrow
  \frac{\ell}{2}^+]{}+\infty 
\end{equation}
Equations \eqref{minor-period}, \eqref{limit-G2G3} and
\eqref{limit-G2t} imply that $\lim F(v_n)=+\infty$. The lemma is
proved. \cqfd

\subsection{Solution of the Period Problem}
In this section, we go back to the $N$-dimensional Period Problem.
Consider $N$ integers $p_1<\cdots<p_N$ (these will specify ``where''
we want to put the handles).
Recall that $\ell<1$ and $\eta=1-\sqrt{1-\ell^2}$ and that lemma
\ref{local-bhv-period} gives us a number $\eta_0<\eta$. Equation
(\ref{eq-condition}) means that each singular point $q_i$ must be at
distance less than $\eta$ from an even integer. We require that
$|q_i-2p_i|<\eta$ for $i=1,\cdots,N$. The next proposition solves the
Period Problem with this setting. 
\begin{proposition}
\label{prop-period-problem}
For any $N$ and for any choice of $p_1,\cdots,p_N$ satisfying
$p_1<p_2<\cdots<p_N$, there exists $(q_1,\cdots,q_N)$ satisfying
$|q_i-2p_i|<\eta_0$ and $F_i(q_1,\cdots,q_N)=0$ for
$i=1,\cdots,N$. 
\end{proposition} 
Proof :
Consider any $p_1,\cdots,p_N$ such that $p_1<\cdots<p_N$.
Consider any value of the parameters $(q_1,\cdots,q_N)$
in the cubic box defined by $q_i\in[2p_i-\eta_0,2p_i+\eta_0]$.
Consider some $i$, $1\leq i\leq N$. Translating by $-2p_i$,
proposition \ref{local-bhv-period} tells us that if $q_i=2p_i+\eta_0$,
$F_i(q_1,\cdots,q_N)\geq 1$ while if $q_i=2p_i-\eta_0$,
$F_i(q_1,\dots,q_N)\leq -1$.  The result then follows from the
Poincar\'e Miranda theorem since $F_i$ depends continuously in the $q_j$.
\cqfd
\begin{remark}
We do not know if the solution to the Period Problem is unique.
Since we do not know how to compute derivatives of the periods with
respect to the parameters, its seems hard to obtain uniqueness.
\end{remark}
\section{Finite genus}
Proposition \ref{prop-period-problem} implies the following
\begin{corollary}
\label{cor-finite-genus}
For each $N\geq 1$, there exists a complete, properly embedded minimal surface
in $\R^2\times\S^1$ which has genus $N$, infinite total curvature,
infinitely many ends, and two limit ends.
\end{corollary}
Proof :
consider integers $p_1<\cdots<p_N$. Let $(q_1,\cdots,q_N)$ be the
solution to the Period Problem given by proposition
\ref{prop-period-problem}. Let $v=v[q_1,\cdots,q_N]$,
$u=u[q_1,\cdots,q_N]$ and $X^*=X^*[q_1,\cdots,q_N]$.
Then $X^*$ is well defined in $\Omega\setminus\{q_1,\cdots,q_N\}$.
To see that the image of $X^*$ is embedded we argue as follows.
Let $M$ be the graph of $u$ on the strip $\Omega^+=\R\times(0,\ell)$
and $M^*$ be the conjugate minimal surface to $M$,
so $M^*=X^*(\Omega^+)$.
Since $\Omega^+$ is convex, $M^*$ is a graph over a planar domain
by the Theorem of R. Krust, so is embedded.

\medskip

Since the Period Problem is solved, all segments $(q_i,q_{i+1})$,
$i=1,\cdots,N-1$, 
as well as the half lines $(-\infty,q_1)$ and $(q_N,+\infty)$,
are mapped onto geodesics in the vertical plane $x=0$ (after a
suitable translation).
Consider now some $(x,y)\in\Omega^+$ such that $x\neq q_i$ for all $i$.
From the formula (\ref{eq-dX1star}) we see that 
$\int_{(x,0)}^{(x,y)} dX_1^*>0$.
Hence the point $X^*(x,y)$ lies in $x>0$. By continuity this is true for
all points of $\Omega^+$, so $M^*$ is an embedded minimal surface in
$x>0$, $0<z<1$. Extending by symmetry with respect to the plane $x=0$
and the horizontal planes at integer heights we obtain an embedded,
periodic minimal surface in $\R^3$ with period $T=(0,0,2)$.
The quotient by this period is an
embedded minimal surface ${\cal M}_N$ in $\R^2\times\S^1$ of
genus $N$ (now and in the following, we identify $\S^1$ with
$\R/2\Z$). We will see in section \ref{section-bounded-curvature} that
it has bounded curvature. By a theorem of Meeks - Rosenberg, theorem
2.1 in \cite{meeksRosenberg},
a complete embedded minimal surface in $\R^3$ with bounded curvature
is properly embedded. Hence ${\cal M}_N$ is properly embedded.
\cqfd

In the following, when $(q_1,\cdots,q_N)$ is given by proposition
\ref{prop-period-problem}, the associated minimal surface given by the
above corollary will be denoted by $\boM[q_1,\dots,q_N]$. This surface
is normalized so that the conjugate to the point $(-1,0,u(-1,0))$
is the point $(0,0,v(-1,0))$.

\begin{figure}[h]
\begin{center}
\resizebox{0.9\linewidth}{!}{\input{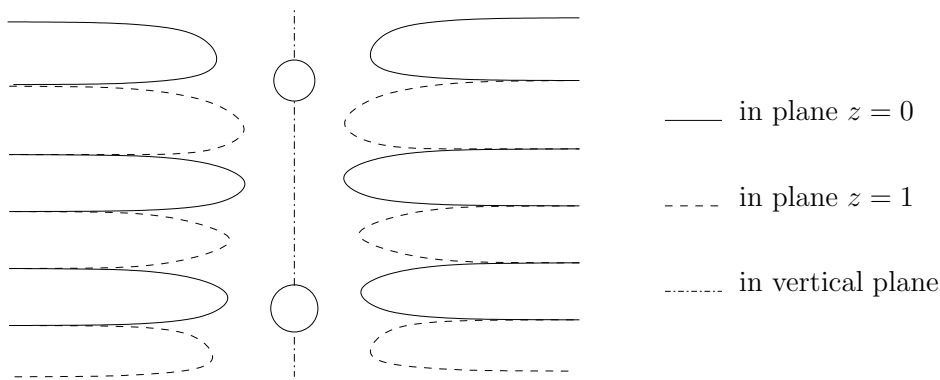}}
\caption{the conjugate surface, Period Problem solved}
\label{fig5}
\end{center}
\end{figure}

\begin{remark}
In proposition \ref{prop-period-problem}, if the singularity set $\X$
is empty, the period problem is solved. Then a surface
$\boM[\emptyset]$ of genus zero exists; in fact this surface is a
Karcher's toroidal halfplane layer.
\end{remark}

\section{Infinite genus}
In this section, we consider the case where we have an infinite
number of singularities.
\begin{proposition}
\label{prop-infinite-genus}
Let $\ell$ and $\eta_0$ be as in proposition
\ref{prop-period-problem}. Consider a strictly increasing sequence of
integers $(p_i)_{i\in\Z}$. Then there exists a sequence
$(q_i)_{i\in\Z}$, such that $|q_i-2p_i|\leq\eta_0$, which solves the
Period Problem $F_j(q_i : i\in\Z)=0$ for all $j\in\Z$.
\end{proposition}
Proof : consider some $n\in\N$. 
We apply proposition \ref{prop-period-problem} with $N=2n+1$ and
the $N$ integers $p_{-n},\cdots,p_n$, and we obtain
$N$ real numbers $q_{-n,n},\cdots,q_{n,n}$ such that
$|q_{i,n}-2p_i|\leq\eta_0$ for all $|i|\leq n$,  
and $F_j(q_{i,n} : |i|\leq n)=0$ for all $|j|\leq n$. 
Using a diagonal process, we can find a subsequence (which we will still
denote by $n$) such that
$\forall i\in\Z$, $q_i=\lim_{n\to\infty} q_{i,n}$ exists.
The limit domain of $(\Omega\setminus\{q_{i,n}:|i|\leq n\})_{n\in\N}$ is
$\Omega\setminus\{q_i:i\in\Z\}$.
Let $v_n=v[q_{i,n} : |i|\leq n]$ and $v_{\infty}=v[q_i:i\in\Z]$.
Arguing as usual, the sequence $(v_n)_{n\in\N}$ has no divergence
line, so up to a subsequence,
it converges on compact subsets of $\overline{\Omega}$
to a solution $v$ of the Dirichlet problem on
$\Omega\setminus\{q_i:i\in\Z\}$.
By uniqueness, $v=v_{\infty}$.
Then $dX_1^*[q_{i,n}:|i|\leq n]$ converges to $dX_1^*[q_i:i\in\Z]$
on compact subsets of $\Omega\setminus\{q_i:i\in\Z\}$.
Integrating on $\gamma_j$ gives
$$F_j(q_i:i\in\Z)=\lim_{n\to\infty} F_j(q_{i,n}:|i|\leq n)=0$$
\cqfd

In the following when a sequence $(p_i)_{i\in\Z}$ satisfies the
hypotheses of the above proposition and $(q_i)_{i\in\Z}$ is a sequence
such that $|q_i-2p_i|\le\eta_0$ for all $i\in\Z$ and
$F_j(q_i:i\in\Z)=0$ for all $j\in\Z$, we shall say that
$(q_i)_{i\in\Z}$ solves the Period Problem for the data
$(p_i)_{i\in\Z}$.

\begin{corollary}
\label{cor-infinite-genus}
For any strictly increasing sequence of integers $(p_i)_{i\in\Z}$,
there exists a properly embedded minimal surface ${\cal M}$ in
$\R^2\times\S^1$ which has infinite genus, infinite total curvature,
infinitely many ends, and two limit ends. Moreover, if the sequence
$(p_{i+1}-p_i)_{i\in\Z}$ is not periodic, then ${\cal M}$ is not
periodic. 
\end{corollary}
Proof : same as proof of corollary \ref{cor-finite-genus}.
\cqfd

As above when $(q_i)_{i\in\Z}$ solves the Period Problem for the data
$(p_i)_{i\in\Z}$, the associated surface is denoted by
$\boM[q_i:i\in\Z]$ and is normalized as in the finite genus case.

Using the notations of proposition \ref{prop-infinite-genus} proof, we
define $\boM_N=\boM[q_{i,N}:-N\le i\le N]$ and
$\boM=\boM[q_i:i\in\Z]$.

\begin{proposition}\label{conv-finite-infinite}
A subsequence of $({\cal M}_N)_N$ converges smoothly 
on compact subsets of $\R^2\times\S^1$ to ${\cal M}$.
\end{proposition}
Proof :
since a subsequence of $(v_n)_{n}$ converges to $v$,
the result seems to be obvious.
This is not immediate for the following reason : the convergence
of the conjugate functions $(u_n)_n$ to $u$ only holds on compact
subsets of $\Omega\setminus\{q_i :i\in\Z\}$. In particular,
this convergence does not tell us anything for the graph of $u$
above the singular points and the vertices $a_k^+$, $k\in\Z$.
Since these correspond to the horizontal symmetry curves on the
conjugate minimal surface, we see that the convergence of
$(v_n)_n$ to $v$ is not enough to conclude. 

\medskip

One way around this difficulty is as follows :
we will see in the next section that the curvature of
${\cal M}_N$ is bounded by a constant independant of $N$.
By the Regular Neighborhood Theorem, or ``Rolling Lemma''
(firstly proven by A.~Ros~\cite{ros5}, Lemma 4, for properly embedded
minimal surfaces in $\R^3$ with  finite total curvature, 
and generalized to properly embedded minimal surfaces with bounded
curvature by Meeks and Rosenberg~\cite{meeksRosenberg}, Theorem 5.3),
each ${\cal M}_N$ has an embedded tubular neighborhood of radius
$1/\sqrt{c}$. In particular, we have local area bounds, namely the
area of ${\cal M}_N$ inside a ball of radius $1/\sqrt{c}$ is bounded
by some constant. By standard results, a subsequence of $({\cal M}_N)_N$
converges smoothly on compact subsets of $\R^2\times\S^1$ to a limit minimal
surface ${\cal M}_{\infty}$. Since $(v_n)_n$ converges to $v$,
the limit ${\cal M}_{\infty}$ needs to be ${\cal M}$.
\cqfd

\medskip

Here is another result in the same spirit, which will be
usefull in section \ref{section-quasi}. For every $n\in\N$,
let $(p_{i,n})_{i\in\Z}$ be a sequence as in Proposition
\ref{prop-infinite-genus}, namely $(p_{i,n})_{i\in\Z}$ is a strictly
increasing sequence of integers. Let $(q_{i,n})_{i\in\Z}$ be a
sequence that solves the Period Problem for the data
$(p_{i,n})_{i\in\Z}$. 
\begin{proposition}\label{conv-inf-inf}
Let  $(p_{i,n})_{i\in\Z}$ and $(q_{i,n})_{i\in\Z}$ be defined as
above. Let us assume that for every $i\in\Z$, $\lim_n
p_{i,n}=p_{i,\infty}$ and $\lim_n q_{i,n}=q_{i,\infty}$. Then
$(q_{i,\infty})_{i\in\Z}$ solves the Period Problem for the data
$(p_{i,\infty})_{i\in\Z}$ and $(\boM[q_{i,n}:i\in\Z])_{n\in\N}$
converges smoothly on compact subsets of $\R^2\times\S^1$ to
$\boM[q_{i,\infty}:i\in\Z]$.  
\end{proposition}

Proof : Let $v_n$ be $v[q_{i,n}:i\in\Z]$. First we notice that the
convergence of $(q_{i,n})_{n\in\Z}$ implies the convergence of
$(p_{i,n})_n\in\Z$. Since, for 
all $i$ and $n$, we have $p_{i+1,n}-p_{i,n}\ge 1$ and
$|q_{i,n}-2p_{i,n}|\le \eta_0$, we obtain $p_{i+1,\infty}-p_{i,\infty}\ge
1$ and $|q_{i,\infty}-2p_{i,\infty}|\le \eta_0$. 

Since all the $(q_{i,n})_{n\in\N}$ converge, the limit domain of
$(\Ome\backslash\{q_{i,n}:i\in\Z\})_{n\in\N}$ is
$\Ome\backslash\{q_{i,\infty}:i\in\Z\}$. As in the preceding, the
sequences $(v_n)_{n\in\N}$ has no divergence line and converges to
a solution
$v_\infty$
of the Dirichlet problem on the limit domain
$\Ome\backslash\{q_{i,\infty}:i\in\Z\}$. By uniqueness,
$v_\infty=v[q_{i,\infty}:i\in\Z]$. As in the proof of Proposition
\ref{prop-infinite-genus},
$$
F_j(q_{i,\infty}:i\in\Z)=\lim_{n\rightarrow\infty}F_j(q_{i,n}:i\in\Z)=0
$$ 
Then $(q_{i,\infty})_{i\in\Z}$ solves the Period Problem. 

Now as in the proof of Proposition \ref{conv-finite-infinite}, since
the curvature of the surfaces $\boM[q_{i,n}:i\in\Z]$ is uniformly
bounded (see Proposition \ref{prop-curvature}), the sequence
$(\boM[q_{i,n}:i\in\Z])_{n\in\N}$ converges 
smoothly on compact subsets of $\R^2\times\S^1$ to a limit minimal
surface $\boM_\infty$. Since $(v_n)_{n\in\N}$ converges to
$v[q_{i,\infty}:i\in\Z]$, the surfaces $\boM_\infty$ needs to be
$\boM[q_{i,\infty}:i\in\Z]$. \cqfd

\section{Bounded curvature}
\label{section-bounded-curvature}
In this section we prove that the curvature of the surfaces $\boM$
given by corollaries \ref{cor-finite-genus} and
\ref{cor-infinite-genus} are bounded by a constant
$C$ depending only on $\ell$. Actually, because of proposition
\ref{conv-finite-infinite}, it suffices to prove it in the finite
genus case. 

\subsection{Size of the handles}
Let $v$ be the solution of the Dirichlet problem on
$\Omega\setminus\{q_1,\cdots,q_N\}$. 
The multi-valuation of its conjugate function $u$ around the
singularity $q_i$ is $\int_{\gamma_i} du$. 
This is equal to twice the length of the vertical segment $B_i$,
which is equal to the length of the geodesic $B_i^*$.
So the multi-valuation of $u$ may be understood as the ``size'' of the
handle. In this section we give a uniform lower bound for this
multi-valuation, which prevents the handles from getting too small.

\medskip

We use the notation $\Omega_L=(-L,L)\times (-\ell,\ell)$.
\begin{proposition}
\label{prop-jump}
Consider some $\ell<1$ and some $\eta_0<\eta$.
There exists $\kappa>0$ (depending on $\ell$ and $\eta_0$) such
that the following is true :
Let $q\in (-\eta_0,\eta_0)$ and $\boS\subset(-2,-2+\eta)\cup(2-\eta,2)$.
Let $v$ be a solution of
the maximal graph equation (\ref{eq-max}) in $\Omega_2\setminus(\{q\}\cup\boS)$
with boundary value
$\varphi$ on $[-2,2]\times\{-\ell,\ell\}$ and $0$ at $\{q\}\cup\boS$.
As in lemma \ref{convergence} the boundary value on the
vertical edges is free, but we require $v$ to be between $0$ and $1$.
Let $u$ be the conjugate function of $v$.
Let $\gamma$ be a small circle around $q$.
Then $\left|\int_{\gamma}du\right|\geq \kappa$.
\end{proposition}
Proof : assume by contradiction that the proposition is not true.
Then there exists sequences $(q_{n})_n$ and $(\boS_n)_n$ and a
sequence $(v_n)_n$ such that $\int_{\gamma}du_n\to 0$. Passing to a
subsequence, $q_{n}$ converges to some $q\in
[-\eta_0,\eta_0]\subset(-\eta,\eta)$. By lemma \ref{convergence}, the
sequence $(v_n)_n$, restricted to $\Omega_1\setminus\{q\}$, does not
have any divergence line, so passing to a subsequence, it converges to
a solution $v$ on $\Omega_1$. 
Then the conjugate differential $du_n$ of $v_n$ converges on
compact subsets of $\Omega_1\setminus\{q\}$ to the conjugate
differential $du$ of $v$. 
This implies that $\int_{\gamma}du=0$, so $u$ is in fact well defined
in $\Omega_1\setminus\{q\}$. Since it satisfies the minimal graph
equation, the 
point $q$ is a removable singularity, so $u$ extends smoothly to
$q$. But then $v$ itself also extends smoothly to $q$. Since
$v(q)=0$ and $0\leq v\leq 1$, the maximum principle for maximal
surfaces gives us that $v=0$ in $\Ome_1$; this contradicts $v=\varphi$
on the boundary.
\cqfd
\subsection{Gradient estimates}
Recall that the graph of $u$ is bounded by a vertical segment
above each singularity $q_i$.
Along this segment, the normal is horizontal.
The following lemma ensures that the normal remains
close to the horizontal on the disk $D(q_i,\delta)$,
where $\delta$ is a number we can control in function of the length of the
vertical segment.
\begin{lemma}
\label{lemma-gradient}
For any $C>0$, for any $\kappa>0$, there exists $\delta>0$ such
that the following is true :
let $v$ be any solution of the maximal graph equation (\ref{eq-max})
on the punctured disk $D(0,1)\setminus\{0\}$ with a singularity at
the origin.
Assume that $v(0)=0$ and $0\leq v\leq 1$.
Let $du$ be the conjugate differential of $v$.
Assume that $|\int_{\gamma}du|\geq \kappa$.
Then $|\nabla u|\geq C$ in $D(0,\delta)\setminus\{0\}$.
\end{lemma}
Proof : assume by contradiction that the lemma is false. Then there exists
$C>0$, $\kappa>0$, and sequences $(v_n)_n$, $(p_n)_n$, such that
$v_n$ is a solution of (\ref{eq-max}) in $D(0,1)\setminus\{0\}$,
$p_n\to 0$, $\int_{\gamma} du_n\geq \kappa$
and $|\nabla u_n(p_n)|\leq C$.
Let $\lambda_n=|p_n|$.
Let $\widetilde{v}_n(p)=v_n(\lambda_n p)/\lambda_n$ (so the
graph of $\widetilde{v}_n$ is the graph of $v_n$ scaled by $1/\lambda_n$).
Let $\widetilde{p}_n=p_n/\lambda_n$; by rotation we may assume that
$\widetilde{p}_n=(1,0)$. 
Then $\widetilde{v}_n$ solves (\ref{eq-max}) in the punctured disk of
radius $1/\lambda_n$. This domain converges to the plane
punctured at the origin.

\medskip

Let us study the convergence of the sequence $(\widetilde{v}_n)_n$.
If there are no divergence lines, then the sequence $(\widetilde{v}_n)_n$
converges on compact subsets of 
the punctured plane to a solution $\widetilde{v}$.
Then the conjugate differentials $d\widetilde{u}_n$ converge to
$d\widetilde{u}=d\Phi_{\widetilde{v}}$. However, 
$$\left|\int_{\gamma}d\widetilde{u}_n\right|
=\frac{1}{\lambda_n}\left|\int_{\gamma}du_n\right|
\geq \frac{\kappa}{\lambda_n}\to\infty$$
so $\int_{\gamma}d\widetilde{u}=\infty$, which is absurd.
So there must be divergence lines.

\medskip

Observe that
$$|\nabla \widetilde{v}_n(1,0)|=|\nabla \widetilde{v}_n(\widetilde{p}_n)|
=|\nabla v_n(p_n)|
=\frac{|\nabla u_n(p_n)|}{\sqrt{1+|\nabla u_n(p_n)|^2}}
\leq \frac{C}{\sqrt{1+C^2}} < 1.$$
Hence the point $(1,0)$ is in the convergence domain of the sequence
$\widetilde{v}_n$.
Let $U$ be the component of the convergence domain containing the
point $(1,0)$. 
Since $\widetilde{v}_n\geq 0$, a divergence line cannot extend infinitely
in both directions, so must be a half-line ending at the origin.
If there are at least two divergence lines then $U$ is a sector defined
in polar coordinates by $0<r<\infty$, $\alpha_1<\theta<\alpha_2$.
The conjugate functions $\widetilde{u}_n$ are well defined in $U$ and converge
to $\widetilde{u}$. Then $\widetilde{u}$ takes the values $\pm \infty$
on the half-lines $\theta=\alpha_1$ and $\theta=\alpha_2$. Since
$\widetilde{v}(0)=0$ and $\widetilde{v}\geq 0$, $\widetilde{u}$
takes the values $+\infty$ on $\theta=\alpha_1$ and $-\infty$ on
$\theta=\alpha_2$. It is proven in \cite{mazet4}, proposition 2, that
this Jenkins Serrin problem has no solution.

\medskip

If there is only one divergence line, then $U$ is a sector of angle $2\pi$
defined in polar coordinates by $0<r<\infty$, $\alpha<\theta<\alpha+2\pi$.
Then $\widetilde{u}$ solves the following Jenkins Serrin problem :
$u=+\infty$ on the half-line $\theta=\alpha$ (approaching this line
with $\theta>\alpha$) and $u=-\infty$ on $\theta=\alpha+2\pi$ (approaching
with $\theta<\alpha+2\pi$.
It is proven in \cite{mazet4}, proposition 4, that this Jenkins Serrin
problem has no solution.
This contradiction proves the lemma.
\cqfd

The following lemma provides a similar estimate in 
a neighborhood of the boundary points $a_k^+$, $k\in\Z$.
It is proven in \cite{mazet-rodriguez-traizet}, lemma 6.
\begin{lemma}
\label{lemma-gradient2}
Given $C>0$, there exists $\delta>0$ such that the following is true :
let $u$ be a solution of the minimal graph equation (\ref{eq-min})
in the half disk $D(0,1)\cap\{y> 0\}$, with boundary values
$+\infty$ on the segment $(0,1)\times\{0\}$ and $-\infty$ on
the segment $(-1,0)\times\{0\}$.
Let $v$ be the conjugate function of $u$. Assume that $v(0)=0$
and $v\geq 0$. Then $|\nabla u|\geq C$ in $D(0,\delta)\cap\{y>0\}$.
\end{lemma}
\subsection{Curvature estimate}

Let $\boM$ be a minimal surface given by corollary
\ref{cor-finite-genus}. 

\begin{proposition}
\label{prop-curvature}
There exists a constant $C$ (depending only on $\ell$)
such that the Gauss curvature $K$ of the surface ${\cal M}$ is bounded by $C$.
\end{proposition}
Proof : since $M^*$ and $M$ are locally isometric, it suffices to bound
the curvature of $M$. The proof is based on stability arguments.
In what follows, all constants involved only depend on $\ell$.

\medskip

Recall that $M$ is the graph of $u$ on
$\Omega^+=\R\times(0,\ell)$.
Let $q_1,\cdots,q_N$ be given by proposition \ref{prop-period-problem}.
By proposition \ref{prop-jump}, there exists $\kappa$ such that
$|\int_{\gamma_i}du|\geq \kappa$ for $i=1,\cdots,N$.
We apply lemma \ref{lemma-gradient} with $C=100$ and obtain
a $\delta_1<\ell$ such that $|\nabla u|\geq 100$ in $D(q_i,\delta_1)$,
$i=1,\cdots,N$.
We apply lemma \ref{lemma-gradient2} with again $C=100$ and obtain
a $\delta_2<\ell$ such that $|\nabla u|\geq 100$ in $D(a^+_k,\delta_2)$,
$k\in\Z$.
We take $\delta=\min\{\delta_1,\delta_2\}$.
Fix some $i=1,\cdots,N$. Let $U$ be the graph of $u$ above the
half disk $D(q_i,\delta)\cap\Omega^+$. Since $|\nabla u|\geq 100$,
the Gauss image of $U$ is included in the spherical domain
$\S^2\cap\{|z|\leq 1/100\}$. The boundary of $U$ consists of a vertical
segment, two horizontal segments and a helix-like looking curve which
is a graph on $\S^1(q_i,\delta)\cap\Omega^+$. Completing by all symmetries,
we obtain a minimal surface $\Sigma$ which is bounded by two helix-like
looking curves, and which is complete in the cylinder $D(q_i,\delta)\times\R$.
The surface $\Sigma$ is of course not a graph anymore. However
its Gauss image is still included in $\S^2\cap\{|z|<1/100\}$.
As the spherical area of this domain is less than $2\pi$, $\Sigma$ is
stable by the theorem of Barbosa Do Carmo \cite{bc1}.
Consider now a point $(x,y)\in D(q_i,\delta/2)$ and let
$p=(x,y,u(x,y))$ be the corresponding point on $M$.
Since $p\in\Sigma$ is at distance more than $\delta/2$ from the boundary of
$\Sigma$, the theorem of Schoen \cite{sc3} ensures that the Gauss curvature
at $p$ is bounded by $c/(\delta/2)^2$ for some universal constant $c$.
The same argument gives the same estimate for the Gauss curvature
when $(x,y)\in D(a^+_k,\delta/2)$, $k\in\Z$.

\medskip

Assume now that $(x,y)\in\Omega^+$ is at distance more than $\delta/2$ from
all points $q_i$ and all points $a^+_k$. 
Let again $p=(x,y,u(x,y))$.
If $y>\delta/4$, then the distance of $p$ to the boundary of $M$
is greater than $\delta/4$ (because $u=\pm \infty$ on the top edges).
Since $M$ is a graph, it is stable, so the Gauss curvature at $p$
is bounded by $c/(\delta/4)^2$.

\medskip

It remains to understand the case $0<y<\delta/4$.
There exists $i$ such that $q_i<x<q_{i+1}$ (with the convention that
$q_0=-\infty$ and $q_{N+1}=+\infty$).
Consider the box $(q_i,q_{i+1})\times(-\delta/2,\delta/2)$. As this
is a simply connected domain of $\Omega$, $u$ is well defined on it.
Let $V$ be the graph of $u$ on this box. The distance of $p=(x,y,u(x,y))$
to the boundary of $V$ is greater than $\delta/4$. Since $V$ is stable,
we conclude again that the Gauss curvature at $p$ is bounded by
$c/(\delta/4)^2$.
\cqfd

\medskip

\section{Quasi-periodicity}
\label{section-quasi}
In this section, we prove that if the sequence
$(p_i-p_{i-1})_{i\in\Z}$ is quasi-periodic but not periodic, then
we can find a solution $(q_i)_{i\in\Z}$ to the Period Problem
such that the associated minimal surface $\boM[q_i\,:\, i\in\Z]$
is quasi-periodic but not periodic in $\R^2\times\S^1$.

\subsection{Preliminaries}
First we need to fix some notation. In the following, an element of
$\R^\Z$ will be denoted as a function: $x\in\R^\Z$ denotes the sequence
$(x(i))_{i\in \Z}$. We consider on $\R^\Z$ the
topology of the pointwise convergence \textit{i.e.} the sequence
$(x^n)_{n\in\N}$ converges to $x^\infty$ if for every $i$ : $\lim_n
x^n(i)=x^\infty(i)$.

We notice that, for every $A\in\R_+$, the subset $[-A,A]^\Z\subset\R^\Z$
is compact. Besides, on $[-A,A]^\Z$, the pointwise convergence is
metrizable : if $x,y\in[-A,A]^\Z$, we define a distance by
$\displaystyle d(x,y)=\sum_{i\in\Z} \frac{1}{2^{|i|}}|x(i)-y(i)|$.

Let $\varphi: \N\rightarrow \N$ be a function. In the following we say that
$\varphi$ is an extraction if $\varphi$ is strictly increasing. The group
$\Z$ acts on the set $\R^\Z$ by shift : if $x\in\R^\Z$ and $n\in\Z$,
we denote by $n\cdot x$ the sequence $(x(n+i))_{i\in\Z}$. Then if
$\varphi$ is an extraction 
and $x\in\R^\Z$, we define the sequence $\varphi\cdot x=(\varphi(n)\cdot
x)_{n\in\N}$ in $\R^\Z$. We have the following definitions.

\begin{definition}
Let $x$ be in $\R^\Z$, this sequence is said to be \emph{quasi-periodic}
if there exists an extraction $\varphi$ such that the sequence
$\varphi\cdot x$ converges pointwise to $x$ (namely, for all $i$,
$\lim_n x(i+\varphi(n))=x_i$).

Let $x\in\R^\Z$ be an increasing sequence, we say that $x$ has
\emph{quasi-periodic gaps} if the sequence $(x(i)-x(i-1))_{i\in\Z}$ is
quasi-periodic. 
\end{definition}

Let us give two examples :
\begin{enumerate}
\item let $\alpha$ be an irrational number, let $x(i)=[\alpha i]$
be the integer part of $\alpha i$ and let $g(i)=x(i)-x(i-1)$.
Then the sequence $(g(i))_{i\in\Z}$ is quasi-periodic and is not
periodic. Moreover, for any extraction $\varphi$,
if $\lim_{n\to\infty}\varphi(n)\cdot g$ exists, then it is not periodic.
\item (the counting sequence) consider the infinite word on the
alphabet $\{0,\cdots,9\}$ formed by writing in order all natural integers :
$$0123456789101112131415161718192021\cdots.$$
For $i\geq 1$, let $x(i)$ be the $i$th
digit in this word. For $i\leq 0$, let $x(i)=0$.
The sequence $(x_i)_{i\in\Z}$ is quasi-periodic but not
periodic.
However, if $\varphi$ is an extraction, the limit of
$\varphi(n)\cdot x$ can very well be periodic (in fact it can be any
sequence of integers between $0$ and $9$).
\end{enumerate}
\subsection{Why are we not done yet ?}
Let us assume that the sequence $(p(i))_{i\in\Z}$ has quasi-periodic
gaps, and let $(q_i)_{i\in\Z}$ be a sequence that solves the Period
Problem for the data $(p(i))_{i\in\Z}$.
We expect the surface ${\cal M}[q(i):i\in\Z]$ to be
quasi-periodic, but unfortunately we cannot prove that.
What we can prove is that there exists a sequence of translations
$T_n$ such that $T_n({\cal M}[q(i):i\in\Z])$ converges
to ${\cal M}[q'(i):i\in\Z]$, where $(q'(i))_{i\in\Z}$ is
another solution to the Period Problem for the same data 
$(p(i))_{i\in\Z}$.
Since we do not know whether the Period Problem has a unique solution,
we cannot ensure that $q'(i)=q(i)$.
(If the reader knows that the solution to the Period Problem is
unique, he may omit what follows. He should also inform the
authors).

\medskip

Our strategy is to prove, using Zorn's lemma, that
amongst all the solutions $(q_i)_{i\in\Z}$ to the Period Problem,
at least one of them yields a quasi-periodic minimal surface.

\subsection{Quasi-periodic surfaces}
Let us consider $\ell$ and $\eta_0$ as in Proposition
\ref{prop-period-problem}.  

Let us now explain how we shall construct a quasi-periodic minimal
surface. Let $p=(p(i))_{i\in\Z}$ be a strictly increasing sequence
with quasi-periodic gaps. In the following, we always assume that
$p(0)=0$. The sequence 
$g=(p(i)-p(i-1))_{i\in\Z}$ is quasi-periodic, we then denote by
$\boA$ the non-empty set of all extractions $\varphi:\N\rightarrow\N$
such that $\displaystyle \lim_{n\rightarrow\infty}\varphi(n)\cdot
g=g$. 

Let us fix a sequence $q=(q(i))_{i\in\Z}$ that solves the Period
Problem for the data $(p(i))_{i\in\Z}$. The problem consists in
building from $q$ a sequence $(q'(i))_{i\in\Z}$ that solve the Period
Problem for the data $(p(i))_{i\in\Z}$ and such that
$\boM[q'(i):i\in\Z]$ is quasi-periodic. 

Let us
denote by $r$ the sequence $q-2p$ : $r(i)=q(i)-2p(i)$ for all
$i\in\Z$. Let $\varphi$ be in $\boA$; $\varphi\cdot r$ is a sequence
of elements of $[-\eta_0,\eta_0]^\Z$. This set is compact so there
exists a subsequence of $(\varphi(n)\cdot r)_{n\in\N}$ that converges
in $[-\eta_0,\eta_0]^\Z$. Thus there exists an extraction $\psi
:\N\rightarrow \N$ such that $(\varphi\circ\psi)\cdot r$ converges. We
notice that, since $\varphi\in\boA$, $\varphi\circ\psi\in\boA$. The
following result decribes such a situation.
\begin{proposition}\label{conv-period-prob}
With the above notation, let $\varphi\in\boA$ such that $\displaystyle
\lim_{n\rightarrow\infty}\varphi(n)\cdot r=r'$. Then
$2p+r'=(2p(i)+r'(i))_{i\in\Z}$ solves the Period Problem for the data
$(p(i))_{i\in\Z}$.
\end{proposition}

Proof : For every $n\in\N$, let us define the sequence $q^n$ by
$q^n(i)=q(i+\varphi(n))-2p(\varphi(n))$ for all
$i\in\Z$. We also define $p^n$ by $p^n(i)=p(i+\varphi(n))-p(\varphi(n))$
for all $i\in\Z$.

\begin{claim}\label{claim2}
We have $\lim p^n=p$ and $\lim q^n=2p+r'$.
\end{claim}

Proof : Let us fix $i\in\Z$ then, if $i\ge 1$:
\begin{align*}
p^n(i)=p(i+\varphi(n))-p(\varphi(n))=
\sum_{l=1+\varphi(n)}^{i+\varphi(n)} p(l)-p(l-1)
&=\sum_{l=1+\varphi(n)}^{i+\varphi(n)}g(l)\\
&=\sum_{l=1}^i\varphi(n)\cdot g(l)
\end{align*}

Since $\lim_{n\rightarrow\infty}\varphi(n)\cdot g=g$, the right-hand term
converges to $\displaystyle \sum_{l=1}^i g(l) =p(i)$. When $i<1$ we
have:
$$p^n(i)=p(i+\varphi(n))-p(\varphi(n))= \sum_{l=i+\varphi(n)+1}
^{\varphi(n)}p(l-1)-p(l)= -\sum_{l=i+1}^0\varphi(n)\cdot g(l)$$
The right-hand term converges again to $p(i)$. Then $\lim p^n=p$.

We have $q^n(i)= r(i+\varphi(n))+2p(i+\varphi(n))- 2p(\varphi(n))=
\varphi(n)\cdot r(i)+2p^n(i)$ for all $i\in\Z$. Since
$\lim_{n\rightarrow\infty}\varphi(n)\cdot r=r'$, $\lim q^n=2p+r$.
\cqfd

\medskip

By definition of $q^n$, the uniqueness of the solution to the Dirichlet
problem implies that we have 
\begin{equation}\label{equa-translat}
v[q^n(i):i\in\Z](x,y)=v[q(i):i\in\Z](x+2p(\varphi(n)),y)
\end{equation}
This implies that $(q^n(i))_{i\in\Z}$ solves the Period Problem :
$F_j(q^n(i):i\in\Z)=0$ for every $j\in\Z$.

Then by Proposition \ref{conv-inf-inf} and Claim \ref{claim2}, $2p+r'$
solves the Period problem for the data $(p(i))_{i\in\Z}$. \cqfd

\medskip

Proposition \ref{prop-infinite-genus} does not give us the uniqueness
of the sequence $(q(i))_{i\in\Z}$ that solves the Period Problem for
the data $(p(i))_{i\in\Z}$, so, as we said in the preceding subsection,
we cannot ensure that the sequences $r$ and $r'$ are equal. Such an
affirmation would be interesting because of the following proposition.

\begin{proposition}\label{quasi-periodic}
With the above notations, if there exists $\varphi\in\boA$ such that
$\displaystyle \lim_{n\rightarrow\infty}\varphi(n)\cdot r=r$, the
surface $\boM[q(i):i\in\Z]$ is quasi-periodic.
\end{proposition}

Proof : We use the notations of the proof of Proposition
\ref{conv-period-prob}. We have the sequences $q^n$, $p^n$. Now Claim
\ref{claim2} says us that $\lim q^n=q$. Let us recall that when
$(a(i))_{i\in\Z}$ solves the Period problem, the surface
$\boM[a(i):i\in\Z]$ is normalized such that the conjugate to the point
in the graph above $(-1,0)$ is the point $(0,0,v(-1,0))$ where
$v=v[a(i):i\in\Z]$.  

As above, \eqref{equa-translat} is true. So our normalization for the
surfaces $\boM$ implies that $\boM[q^n(i):i\in\Z]$ is the image of
$\boM[q(i):i\in\Z]$ by an horizontal tranlation $T_n$. The vector of
the translation is $(0,-X_2^*(2p(\varphi(n))-1))$ where $X_2^*$ is
$X_2^*[q(i):i\in\Z]$.  

Then by Proposition \ref{conv-inf-inf} and Claim \ref{claim2}, the
sequence of minimal surfaces $(\boM[q^n(i):i\in\Z])_{n\in\N}=
(T_n(\boM[q(i):i\in\Z]))_{n\in\N}$ converges to $\boM[q(i):i\in\Z]$
smoothly on compact subsets of $\R^2\times\S^1$. Since
$\boM[q(i):i\in\Z]$ is properly embedded, $\lim
|X_2^*(2p(\varphi(n))-1)|=\infty$ ; thus $(T_n)_{n\in\N}$ is a
diverging sequence of translations. This proves that
$\boM[q(i):i\in\Z]$ is quasi-periodic. \cqfd

\medskip

By using a proposition that will be proved in the next subsection we
then can prove our main theorem.

\begin{theorem}\label{quasi-periodic-surface}
Let $(p(i))_{i\in\Z}$ be a sequence with quasi-periodic gaps that
satisfies the hypotheses of Proposition
\ref{prop-infinite-genus}. Then there exists a sequence
$(q(i))_{i\in\Z}$ which solves the Period Problem for the data
$(p(i))_{i\in\Z}$ and such that $\boM[q(i):i\in\Z]$ is
quasi-periodic. Moreover if $(p(i+1)-p(i))_{i\in\Z}$ is not periodic,
$\boM[q(i):i\in\Z]$ is not periodic.
\end{theorem}

Proof : By Proposition \ref{prop-infinite-genus}, there exists a
sequence $(q(i))_{i\in\Z}$ that solves the Period Problem for the data
$(p(i))_{i\in\Z}$. 

The sequence $(g(i))_{i\in\Z}=(p(i)-p(i-1))_{i\in\Z}$ is
quasi-periodic so we have the set $\boA$. Let $r$ denotes the sequence
$q-2p$, we recall that $r\in[-\eta_0,\eta_0]^\Z$. By Proposistion
\ref{dynamic}, there exists $\varphi$ and $\psi\in \boA$ such that
\begin{align}
&\lim \varphi\cdot r=r'\label{equ1}\\
&\lim \psi\cdot r'=r'\label{equ2}
\end{align}

By Proposition \ref{conv-period-prob}, equation \eqref{equ1} implies that the
sequence $2p+r'$ solves the Period Problem for the data
$(p(i))_{i\in\Z}$. Equation \eqref{equ2} gives us by Proposition
\ref{quasi-periodic} that $\boM[2p(i)+r'(i):i\in\Z]$ is
quasi-periodic. \cqfd

\subsection{A dynamical result}

Let $X$ be a topological space with a countable basis. In the
following, we shall denote by $(V_n(x))_{n\in\N}$ a countable
decreasing basis of open neighborhoods of $x\in X$. Let
$F:X\rightarrow X$ be a continuous map. Let $g$ be in $X$. We assume
that there exists an extraction $\varphi$ such that $\lim_n
F^{\varphi(n)}(g)=g$. As above we denote by $\boA$ the set of
extractions $\varphi$ such that $\lim_n F^{\varphi(n)}(g)=g$. The aim
of this section is to prove the following proposition. 

\begin{proposition}\label{dynamic}
Let $K$ be a compact subset of $X$ such that $F(K)\subset K$. Let $x$
be in $K$. Then there exists two extractions $\varphi\in\boA$ and $\psi\in\boA$
such that $\lim_n F^{\varphi(n)} (x)=x'$ and $\lim_n
F^{\psi(n)}(x')=x'$.  
\end{proposition}

In the proof of Theorem \ref{quasi-periodic-surface}, we use this result with
$X=\R^\Z$ with its pointwise convergence topology, $K$ is
$[-\eta_0,\eta_0]^\Z$, $F$ is the shift map and $g$ is the
quasi-periodic sequence $g$.

Before proving the above proposition, let us fix some notations. Let
$x$ be as in the proposition and $\varphi\in\boA$; the sequence 
$F^{\varphi(n)} (x)$ is a sequence in $K$ which is compact. Thus
there exists a subsequence that converges. As said above this implies
that there exists an extraction $\psi$ such that $F^{\varphi(\psi(n))}
(x)$ converges. We notice that $\varphi\circ \psi\in\boA$. Hence we
define:  
$$Asymp(x)=\{y\in K\,|\, \exists\varphi\in\boA,\ y=\lim_n
F^{\varphi(n)} (x)\}$$  

We know that $Asymp(x)$ is non-empty. In fact Proposition
\ref{dynamic} consists in proving that there exists $x'\in Asymp(x)$
such that $x'\in Asymp(x')$. Then Proposition \ref{dynamic} is a
consequence of the following three lemmae.

\begin{lemma}
Let $x\in K$, $Asymp(x)$ is a closed subset of $K$.
\end{lemma}

Proof : Let $(y_k)_{k\in\N}$ be a sequence in $Asymp(x)$ that
converges to $y\in K$. For each $k$, we choose
$\varphi_k\in\boA$ such that $y_k=\lim_n F^{\varphi_k(n)}(x)$. We are going
to construct by induction $\psi\in\boA$ such that $y=\lim_n
F^{\psi(n)} (x)$. 

Let $n$ be in $\N^*$, we assume that $\psi(q)$ is constructed for
$q<n$ such that for every $q<n$:
\begin{align*}
F^{\psi(q)} (g)&\in V_q(g),&
F^{\psi(q)} (x)&\in V_q(y)
\end{align*}
Since $\lim y_k=y$, there exists $k_0$ such that
$y_{k_0}\in V_n (y)$; hence $V_n(y)$ is an open neighborhood of
$y_{k_0}$. Since $\varphi_{k_0}\in\boA$, there exists $q_0$ such that
$\varphi_{k_0}(q_0)> \psi(n-1)$ and  
\begin{align*}
F^{\varphi_{k_0}(q_0)}(g)&\in V_n(g),&
F^{\varphi_{k_0}(q_0)}(x)&\in V_n(y)
\end{align*}
Then if we take $\psi(n)=\varphi_{k_0}(q_0)$ we get:
\begin{gather}
F^{\psi(n)}(g)\in V_n(g)\label{Equ1}\\
F^{\psi(n)}(x)\in V_n(y)\label{Equ2}
\end{gather}
This finishes our construction.

Equation \eqref{Equ1} implies that $\psi\in\boA$ and \eqref{Equ2} implies
that $\lim_n F^{\psi(n)}(x)=y$ thus $y\in Asymp(x)$. \cqfd

\medskip

\begin{lemma}
Let $x\in K$, let $y$ be in $Asymp(x)$; then $Asymp(y)\subset
Asymp(x)$.  
\end{lemma}

Proof : Let $z$ be in $Asymp(y)$. Let $\varphi$ and $\psi\in\boA$ such
that $\lim_n F^{\varphi(n)}(x)=y$ and $\lim_n F^{\psi(n)}(y)=z$. Let
us build by induction $\chi\in\boA$ such that $\lim_n
F^{\chi(n)}(x)=z$. 

Let $n$ be in $\N^*$, we assume that $\chi(q)$ is constructed for
$q<n$ such that for every $q<n$:
\begin{align*}
F^{\chi(q)}(g)&\in V_q(g),&
F^{\chi(q)}(x)&\in V_q(z)
\end{align*}
Since $\lim_n F^{\psi(n)}(y)=z$, there exists $k_0$ such that:
\begin{align*}
F^{\psi(k_0)}(g)&\in V_n(g),&
F^{\psi(k_0)}(y)&\in V_n(z)
\end{align*}
Then $(F^{\psi(k_0)})^{-1}\big(V_n(g)\big)$ is an open neigborhood of
$g$ and $(F^{\psi(k_0)})^{-1}\big(V_n(z)\big)$ is an open neighborhood
of $y$. Since $\lim_n F^{\varphi(n)} (x)=y$, there exists $l_0$ such that
$\psi(k_0)+\varphi(l_0)>\chi(n-1)$ and 
\begin{align*}
F^{\varphi(l_0)} (g)&\in (F^{\psi(k_0)})^{-1}\big(V_n(g)\big)&
F^{\varphi(l_0)} (x)&\in (F^{\psi(k_0)})^{-1}\big(V_n(z)\big)
\end{align*}
Hence if we take $\chi(n)=\psi(k_0)+\varphi(l_0)$, we have:
\begin{gather}
F^{\chi(n)} (g)\in V_n(g)\label{Equ3}\\
F^{\chi(n)} (x)\in V_n(z)\label{Equ4}
\end{gather}
This finishes our construction.

Equation \eqref{Equ3} implies that $\chi$ is in $\boA$ and \eqref{Equ4} gives
us that $\lim_n F^{\chi(n)} (x)=z$, hence $z\in Asymp(x)$. \cqfd

\medskip

\begin{lemma}\label{zorn}
Let $K$ be a compact set and $T:K\rightarrow {\mathcal{P}}(K)$ a map
such that:
\begin{enumerate}
\item for all $x\in K$, $T(x)$ is closed and non-empty.
\item for all $x\in K$ and all $y\in T(x)$, $T(y)\subset T(x)$.
\end{enumerate}
Let $x\in K$, then there exists $y\in T(x)$ such that $y\in T(y)$.
\end{lemma}

Proposition \ref{dynamic} is then a consequence of this lemma with
$T=Asymp$.

Proof of Lemma \ref{zorn} : The proof of this lemma is given by Zorn's
Lemma. Let $x$ be in $K$, we denote by $\boB$ the set $\{T(y),y\in
T(x)\}$. $\boB$ is ordered by the inclusion. Let $(T_i)_{i\in I}$ be a
totally ordered family of $\boB$. Let us define $T_\infty=\bigcap_{i\in
I}T_i$. If $T_\infty$ is empty, since each $T_i$ is closed and $K$ is
compact there exists a finite subset $I_0\subset I$ such that
$\bigcap_{i\in I_0}T_i=\emptyset$. Since $(T_i)_{i\in I_0}$ is totally
ordered there exists $i_0\in I_0$ such that $T_{i_0}=\bigcap_{i\in
I_0}T_i$, but $T_{i_0}$ is non-empty thus $T_\infty\neq \emptyset$. 

Let $y$ be in $T_\infty$, then $y\in T_i$ for all $i\in I$. This
implies by the second hypothesis that $y\in T(x)$ and
$T(y)\in\boB$. Besides $T(y)\subset T_i$ for all $i\in I$; then
$T(y)\subset T_\infty$. We obtain that $T(y)$ is an under-bound for
the family $(T_i)_{i\in I}$.

We have proved that every totally ordered family admits an
under-bound. Hence, by Zorn's Lemma, there exists an element $T_m\in\boB$
which is minimal for the inclusion. Let $y$ be in $T_m$ (we
recall that all elements of $\boB$ are non-empty subsets of $K$). We
have $y\in T(x)$ by the second hypothesis  then $T(y)\in\boB$ and
$T(y)\subset T_m$. Since $T_m$ is minimal in $\boB$, $T(y)=T_m$ and
$y\in T(y)$. \cqfd

\bibliographystyle{plain}
\bibliography{bill}
\end{document}